\documentclass[conference,a4paper,10pt]{IEEEtran}
\IEEEoverridecommandlockouts
\usepackage{cite}
\usepackage{amsmath}    
\usepackage{amssymb,amsthm}    
\usepackage{algorithmic}
\usepackage{graphicx}
\usepackage{textcomp}
\usepackage{xcolor}

\usepackage{chemformula}
\usepackage{bm}
\usepackage{subcaption}
\usepackage{siunitx}
\usepackage{flushend}

\usepackage{geometry}
\graphicspath{{./figures/}}

\DeclareSIUnit\cell{cell}
\DeclareSIUnit\cells{cells}

\newcommand{\T}{\mathsf{T}}
\renewcommand{\L}{\mathcal{L}}
\renewcommand{\b}{\boldsymbol}

\newcommand{\x}{\b{x}}
\newcommand{\w}{\b{w}}

\def\BibTeX{{\rm B\kern-.05em{\sc i\kern-.025em b}\kern-.08em
		T\kern-.1667em\lower.7ex\hbox{E}\kern-.125emX}}
\begin{document}

\title{A Meshless Solution of a Small-Strain Plasticity Problem\\
	\thanks{The authors would like to acknowledge the financial support of
		Slovenian Research Agency (ARRS) in the framework of the research core funding
		No. P2-0095 and the World Federation of Scientists.}
}

\author{\IEEEauthorblockN{Filip Strni\v{s}a\textsuperscript{a}, Mitja Jan\v{c}i\v{c}\textsuperscript{a,b}, Gregor Kosec\textsuperscript{a}}
	\IEEEauthorblockA{\textsuperscript{a} Jo\v{z}ef Stefan Institute, Parallel and Distributed Systems Laboratory, Ljubljana, Slovenia\\
		\textsuperscript{b} Jo\v{z}ef Stefan International Postgraduate School, Ljubljana, Slovenia\\
		filip.strnisa@ijs.si, mitja.jancic@ijs.si, gregor.kosec@ijs.si}
}

\maketitle

\begin{abstract}
	When the deformations of a solid body are sufficiently large, parts of the body undergo a permanent deformation commonly refereed to as plastic deformation. Several plasticity models describing such phenomenon have been proposed, e.g. von Mises, Tresca, etc. Traditionally, the finite element method (FEM) is the numerical tool of choice for engineers who are solving such problems. In this work, however, we present the implementation of the von Mises plasticity model with non-linear isotropic hardening in our in-house developed MEDUSA library, utilizing a variant of meshless methods -- namely the radial basis function-generated finite differences (RBF-FD). We define a simple plane stress case, where a 2D block is fixed at one edge, and a tensile force, which causes the block to deform, is applied to it at the opposite edge. We show that results are in good agreement with the numerical solution obtained by Abaqus FEA, a commercial FEM solver.
\end{abstract}

\bigskip


\section{Introduction}
Speaking very broadly, a deformation of a solid body can be broken down into two main sub-categories: elastic and plastic deformation. It is said, that the deformation is elastic, if the body returns to its original shape after the applied load had been released, while plastic deformation occurs when any part of a solid body undergoes a non-reversible change of shape due to sufficiently large load applied~\cite{fung2001classical}. Generally, the material response beyond the elastic-plastic tipping point, commonly referred to as \emph{yielding criterion}, is non-linear~\cite{desouza2008}, thus numerical treatment of partial differential equations is required~\cite{desouza2008}.

Traditionally, such problems are solved with the finite elements method (FEM)~\cite{Bartels2012QuasiStaticSP,Schrder2015SmallSP,Roostaei2018ACS,amouzou2021}. In this work, however, we employ meshless methods that have proven to be a good alternative as they can operate on nodes contrary to mesh-based methods that require meshes~\cite{belytschko1996meshless}. An often used variant of the meshless methods is the radial basis function-generated finite differences (RBF-FD)~\cite{bayona2017role}, which has already been employed to obtain solutions to elasticity~\cite{Depolli2019,slak2019rbffd} and plasticity~\cite{JANKOWSKA201812, jiang2021nonlinear} problems. 

We present our implementation of a von Mises plasticity model with non-linear isotropic hardening limited to small strains in a plane stress example. The implementation was done using our in-house developed MEDUSA C++ library~\cite{slak2021} supporting all the required meshless procedures. The original FEM formulation of the solution procedure provided by the de Souza et al.~\cite{desouza2008} is adapted to employ RBF-FD and used to solve a simple plane stress problem.



\section{Plasticity}
Plasticity problems, where an external force acts on a solid body, are usually solved by applying partial loads to the system i.e. the external force is applied incrementally~\cite{desouza2008,yarushina2010}.
At each increment of the external force one first predictably solves the steady-state Navier-Cauchy equation for elastic bodies, and corrects the solution if the plastic yield criterion is violated.
This section will present the small-strain plasticity model for plane stress cases.
The steady state Navier-Cauchy equation can thus be written with Lam\'{e} constants as:

\begin{equation}
	\label{eq: nce}
	\left(\frac{2 \lambda \mu}{\lambda + 2 \mu} + \mu\right) \nabla (\nabla \cdot \bm{u}) + \mu \nabla^2 \bm{u} = -\bm{r},
\end{equation}

\noindent
where $\bm{u} = (u_x, u_y)$ is the displacement vector.
$\lambda = \frac{E \nu}{(1 - 2 \nu) (1 + \nu)}$ and $\mu = \frac{E}{2 (1 + \nu)}$ are the Lam\'{e} parameters, which depend on Young's modulus $E$, and Poisson's ratio $\nu$.
$\bm{r}$ is the residual force density.
The latter is computed as the difference between internal and external force densities: $\bm{r} = \bm{f}^{i} - \bm{f}^e$.
Internal force density is defined as:

\begin{equation}
	\label{eq: sigdiv}
	\bm{f}^i = \nabla \cdot \sigma,
\end{equation}

\noindent
where $\sigma$ is the stress tensor with components: $\sigma_{xx}$, $\sigma_{yy}$, and $\sigma_{xy} = \sigma_{yx}$.
Other components are assumed to equal 0 when assuming plane stress.
For simplicity $\sigma$ can be reshaped into a vector $\sigma = (\sigma_{xx},\sigma_{yy},\sigma_{yx})$.
$\sigma$ is related to $\bm{u}$ \textit{via} the elastic strain tensor $\varepsilon^E$ and the elastic stiffness tensor $\mathrm{D}$ as $\sigma = \mathrm{D} \varepsilon^E$.
$\varepsilon^E$ can be reshaped as $\varepsilon^E = (\varepsilon^E_{xx},\varepsilon^E_{yy},2\varepsilon^E_{yx})$, and then $\mathrm{D}$ becomes:

\begin{equation}
	\label{eq: Dten}
	\mathrm{D} = \left(
	\begin{matrix}
			2 \mu + \frac{2 \lambda \mu}{\lambda + 2 \mu} & \frac{2 \lambda \mu}{\lambda + 2 \mu}         & 0   \\
			\frac{2 \lambda \mu}{\lambda + 2 \mu}         & 2 \mu + \frac{2 \lambda \mu}{\lambda + 2 \mu} & 0   \\
			0                                             & 0                                             & \mu
		\end{matrix}
	\right).
\end{equation}

\noindent
In plasticity $\varepsilon^E$ is a part of total strain tensor $\varepsilon = (\varepsilon_{xx},\varepsilon_{yy},2\varepsilon_{yx}) = \varepsilon^E + \varepsilon^P$, where $\varepsilon^P = (\varepsilon^P_{xx},\varepsilon^P_{yy},2\varepsilon^P_{yx})$ is the irreversible plastic strain.
$\varepsilon$ is computed as:

\begin{equation}
	\label{eq: tot_strain}
	\varepsilon = \frac{\nabla \bm{u} + (\nabla \bm{u})^{\intercal}}{2}.
\end{equation}

As the force is increased incrementally by $\delta \bm{f}^e$, Eq. (\ref{eq: nce}) is used to compute the displacement increment $\delta \bm{u}$.
Thus after $n$ increments $\bm{u}^n = \bm{u}^{n-1} + \delta \bm{u}$.
After finding $\bm{u}^n$ with the elastic prediction, one computes $\sigma$, and checks the yield criterion, in this case the von Mises criterion $\Phi$~\cite{desouza2008}:

\begin{equation}
	\label{eq: vms}
	\Phi = \frac{1}{2} \sigma^{\intercal} \mathrm{P} \sigma - \frac{1}{3} \sigma_Y^2(\varepsilon^P_{eq}) = \frac{1}{2} \xi - \frac{1}{3} \sigma_Y^2(\varepsilon^P_{eq}),
\end{equation}

\noindent
where $\varepsilon^P_{eq}$ is the scalar, ``equivalent'' plastic strain, and $\sigma_Y$ is the likewise scalar yield stress, and is defined by the yield function.
It can be noted, that von Mises stress is $\sigma_{VM} = \sqrt{\frac{3}{2} \sigma^{\intercal} \mathrm{P} \sigma}$.
$\mathrm{P}$ is defined as:

\begin{equation}
	\label{eq: P}
	\mathrm{P} = \frac{1}{3}\left(
	\begin{matrix}
			2  & -1 & 0 \\
			-1 & 2  & 0 \\
			0  & 0  & 6
		\end{matrix}
	\right).
\end{equation}

\noindent
Where the yield criterion is violated ($\Phi > 0$), $\sigma$ is then corrected \textit{via} a local Newton-Rhapson method.
Let there be a correction factor $\delta \gamma$.
First it is set to $\delta \gamma = 0$, then the solver updates it with:

\begin{equation}
	\label{eq: dgam}
	\delta \gamma^{new} = \delta \gamma^{old} - \frac{\Phi}{\Phi'}.
\end{equation}

\noindent
$\Phi' = \frac{1}{2} \xi' - \frac{1}{3} H'$ is the derivative of the yield criterion, where:

\begin{equation}
	\label{eq: xistr}
	\xi' = -\frac{(\sigma_{xx} + \sigma_{yy})^2}{9 \left(1 + \frac{E \delta \gamma}{3(1 - \nu)}\right)^3} \frac{E}{1 - \nu} - 2 \mu \frac{(\sigma_{yy} - \sigma_{xx})^2 + 4 \sigma_{xy}^2}{(1 + 2 \mu \delta \gamma)^3},
\end{equation}

\begin{equation}
	\label{eq: dh}
	H' = 2 \sqrt{\frac{2}{3}} H \left[\sqrt{\xi} + \frac{\delta \gamma \xi'}{2 \sqrt{\xi}}\right] \sigma_Y \left(\varepsilon^P_{eq} + \delta \gamma \sqrt{2 \frac{\xi}{3}}\right).
\end{equation}

\noindent
$H$ is the slope of $\sigma_Y$ at $\varepsilon^P_{eq} + \delta \gamma \sqrt{2 \frac{\xi}{3}}$.
Afterwards $\xi$ is updated:

\begin{equation}
	\label{eq: xiup}
	\xi = \frac{(\sigma_{xx} + \sigma_{yy})^2}{6 \left(1 + \frac{E \delta \gamma}{3 (1 - \nu)}\right)^2} + \frac{(\sigma_{yy} - \sigma_{xx})^2 + 4 \sigma_{xy}^2}{2 (1 + 2 \mu \delta \gamma)^2}.
\end{equation}

\noindent
$\Phi$ (Eq. (\ref{eq: vms})) is also updated, and Eqs. (\ref{eq: dgam})--(\ref{eq: xiup}) are iterated over until the yield criterion is satisfied ($\left|\Phi\right| < \mathrm{tolerance}$).
This is the local iteration.

With computed $\delta \gamma$ the variables are updated as:
\begin{equation}
	\label{eq: update}
	\begin{split}
		\sigma^{new} &= \mathrm{A} \sigma^{old}, \\
		\varepsilon^E &= \mathrm{D}^{-1} \sigma^{new}, \\
		\varepsilon^{P \ new}_{eq} &= \varepsilon^{P \ old}_{eq} + \delta \gamma \sqrt{2 \frac{\xi}{3}}, \\
		\varepsilon^{P \ new} &= \varepsilon^{P \ old} + \delta \gamma \mathrm{P} \sigma^{new}.
	\end{split}
\end{equation}

\noindent
$\mathrm{A}$ is a $\delta \gamma$ dependent tensor, defined as:

\begin{equation}
	\label{eq: A tens}
	\mathrm{A} = \left(
	\begin{matrix}
			\frac{1}{2} (a_1 + a_2) & \frac{1}{2} (a_1 - a_2) & 0   \\
			\frac{1}{2} (a_1 - a_2) & \frac{1}{2} (a_1 + a_2) & 0   \\
			0                       & 0                       & a_2
		\end{matrix}
	\right),
\end{equation}

\noindent
where $a_1 = \frac{1 - \nu}{1 - \nu + \frac{1}{3} E \delta \gamma}$, and $a_2 = \frac{1}{1 + 2 \mu \delta \gamma}$.
With variables updated, $\bm{f}^i$ is updated \textit{via} Eq. (\ref{eq: sigdiv}), and used in Eq. (\ref{eq: nce}) to recompute $\delta \bm{u}$.
This is repeated until $\left|\bm{r}\right| < \mathrm{tolerance}$.
This is the global iteration.
Once the global iteration is completed, $\bm{f}^e$ is increased by $\delta \bm{f}^e$.
The whole process is repeated until $\bm{f}^e$ equals the prescribed value.

Such approach, where $\mathrm{D}$ remains unchanged throughout the computation, requires many global iterations.
This can be improved by exchanging the elastic stiffness tensor with the consistent tangent operator, which is updated by the local iteration~\cite{desouza2008,yarushina2010}.
However, for simplicity this is not explored in this work.

\section{RBF-FD approximation of differential operators}
\label{sec:rbffd}
Often mentioned advantage of the mesh-free methods over the mesh-based methods is that they can operate on scattered nodes. This is particularly convenient when complex three-dimensional domains are being treated, because in such cases, automated mesh generation is even today impossible without a human interference. While there are also other advantages, such as direct control over the approximation order, it is important to know that in general, mesh-free methods are computationally more complex because larger support sizes are needed.

Since the emergence of mesh-free methods, many approximation variants have been proposed. The first mentions of RBF-FD reach back in 2000 with the introduction from Tolstykh~\cite{tolstykh2000using}. Since then, the method has been thoroughly studied and applied to numerous real-life problems. 

In the context of RBF-FD, a linear differential operator $\L$ (for example operators in the Navier-Cacuhy equation~\eqref{eq: nce}) is approximated in node $\x_c$ over set of $n$ nearby nodes
\begin{eqnarray}
	\label{eq:approx}
	\widehat{\L u}(\x_c)=\sum_{i=1}^nw_iu(\x_i)
\end{eqnarray}
for an arbitrary function $u$. A set of $n$ neighboring nodes is often also referred to as the \emph{stencil} nodes or \emph{support} nodes. The unknown weights $\b w$ are obtained for a given set of radial basis functions (RBFs) $\phi$ centered at the stencil nodes of a central node $\x _c$
\begin{eqnarray}
	\phi(\x) = \phi(\left\| \x - \x_c\right \|).
\end{eqnarray}
The approximation~\eqref{eq:approx} can then be written in a linear system
\begin{equation}
	\underbrace{
		\begin{bmatrix}
			\phi_1(\x_1) & \cdots & \phi_n(\x_1) \\
			\vdots       & \ddots & \vdots       \\
			\phi_1(\x_n) & \cdots & \phi_n(\x_n) \\
		\end{bmatrix}
	}_{\b \Phi}
	\underbrace{
		\begin{bmatrix}
			w_1    \\
			\vdots \\
			w_n    \\
		\end{bmatrix}}_{\b w} =
	\underbrace{
		\begin{bmatrix}
			(\L \phi_1(\x)\big|_{ \b x = \b x_c} \\
			\vdots                               \\
			(\L \phi_n(\x)\big|_{ \b x = \b x_c} \\
		\end{bmatrix}}_{\ell_\phi}.
\end{equation}
While researchers in the early work on RBF-FD used infinitely smooth RBFs, such as Gaussians or Multiquadrics, nowadays, piecewise smooth polyharmonic splines (PHS)
\begin{equation}
	\phi(r) = \begin{cases}r^k,       & k \text{ odd}  \\
             r^k\log r, & k \text{ even}\end{cases}
\end{equation}
are commonly used. These, however, do no ensure the convergent behavior nor solvability of the system~\cite{bayona2017role}. That is handled by additionally enforcing the constraint~\eqref{eq:approx} for a set of $s = \binom{m+d}{d}$ monomials with up to and including degree $m$ in a $d$-dimensional domain. 

With the additional constraints, the approximation can be compactly written as
\begin{equation}
	\label{eq:rbf-system-aug}
	\begin{bmatrix}
		{\b \Phi} & {\b P} \\
		{\b P}^\T & \b 0
	\end{bmatrix}
	\begin{bmatrix}
		\b w \\
		\b \lambda
	\end{bmatrix}
	=
	\begin{bmatrix}
		\b\ell_\phi \\
		\b\ell_p
	\end{bmatrix},
\end{equation}
where $\b P$ is an $n\times s$ matrix of monomials evaluated at stencil points, $\b \ell_p$ is the vector of values assembled by applying the considered operator $\L$ to the polynomials at $\x_c$, i.e\ $\ell_p^i = (\L p_i(\x))\big|_{ \b x = \b x_c}$ and $\b \lambda $ are Lagrangian multipliers. Finally, the weights $\w$ to approximate a differential operator $\L$ can be obtained by solving the system~\eqref{eq:rbf-system-aug} using a standard solver from the Eigen library and Lagrangian multipliers are discarded.

\section{Example}
\begin{figure}
	\centering
	\includegraphics[width=\columnwidth]{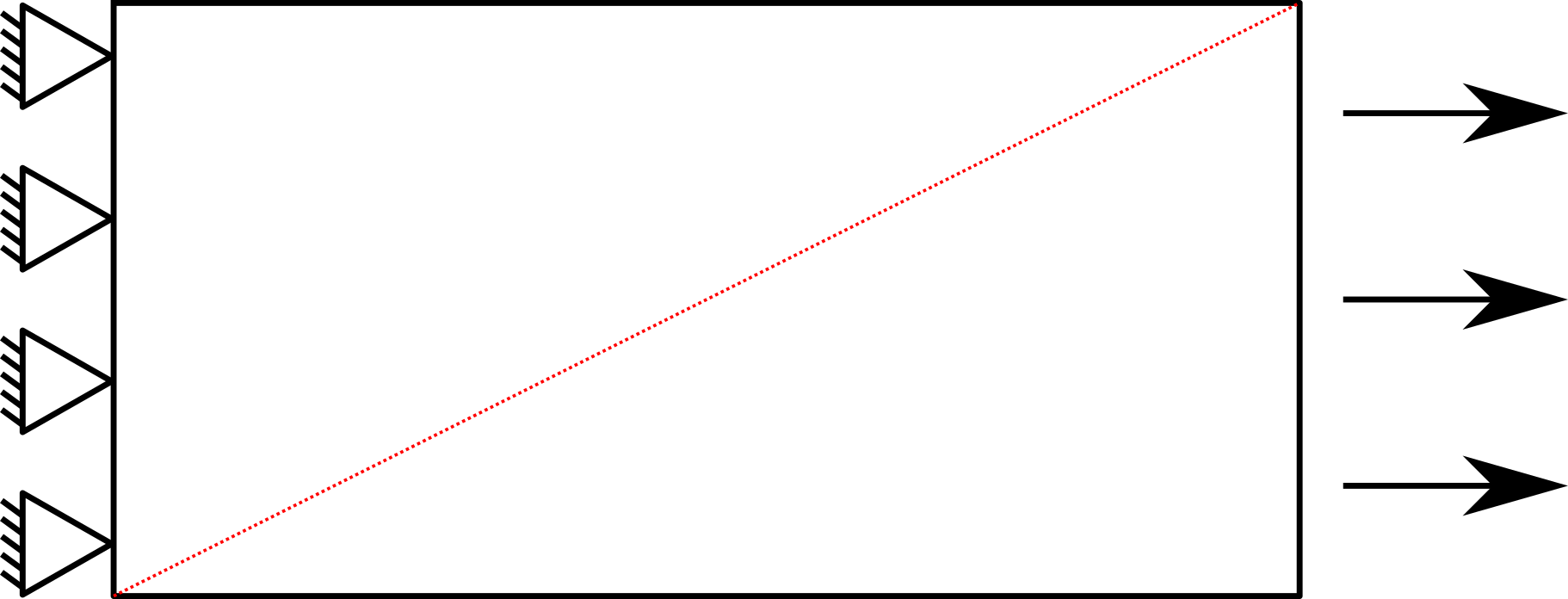}
	\caption{The sketch of the test case. The arrows on the east side represent the tensile stress. The red dashed line represents the results sampling line, used in later analyses. It runs from the south-west corner, to the north-east one.}\label{fig: sketch full}
\end{figure}

For demonstration purposes, a simple plane stress problem is solved. We chose a 2D block (a rectangle), attached to a solid wall on the west edge ($\bm{u} = (0, 0)$) and pulled by a tensile stress on the east edge ($\sigma_{xx} = \SI{30}{\mega\pascal}$), with remaining edges being traction-free ($\bm{n}\cdot\bm{\sigma} = \bm{0}$, where $\bm{n}$ is the face normal), as is shown by the sketch in Fig.~\ref{fig: sketch full}.

The tensile stress acts as the external force on the system, and is increased in steps of equal magnitude, as prescribed by the model input.
Dimensions of the block are as follows: length $L=\SI{10}{\milli\meter}$, and height $H=\SI{5}{\milli\meter}$.
Material properties in this case are chosen arbitrarily, and are: $E = \SI{10}{\giga\pascal}$, $\nu = 0.4$, and a yield function $\sigma_Y(\varepsilon^P_{eq})$ defined by the following set of points ($\varepsilon^P_i, \sigma_{Y i}$): (0, \SI{20}{\mega\pascal}), (0.001, \SI{25}{\mega\pascal}), (0.005, \SI{30}{\mega\pascal}), (0.02, \SI{40}{\mega\pascal}) with a piecewise linear interpolation in-between.

\begin{figure}
	\centering
	\includegraphics[width=\columnwidth]{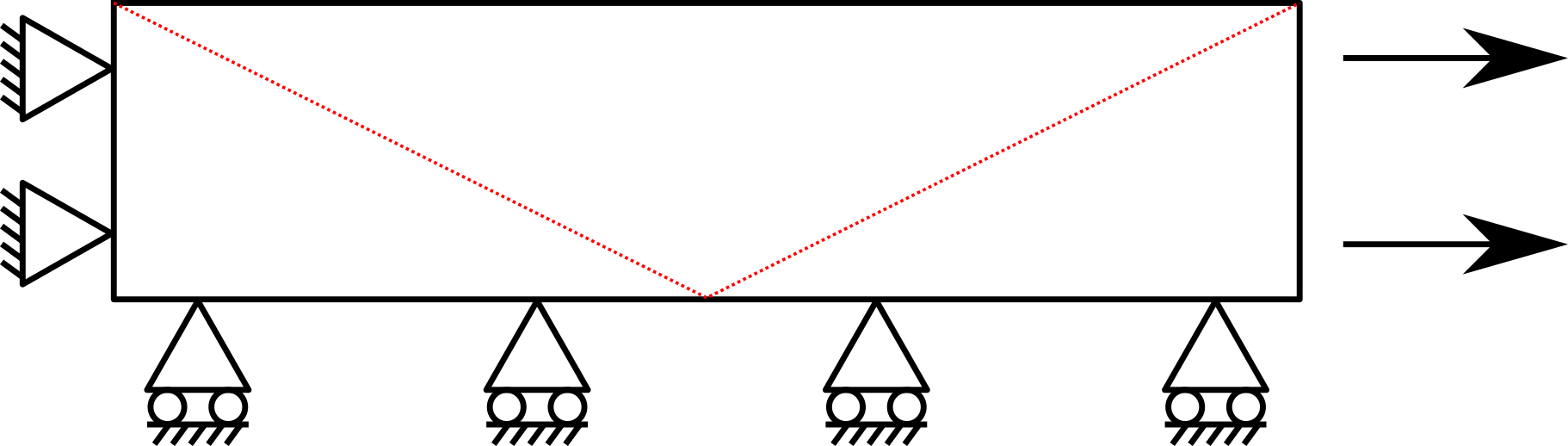}
	\caption{The sketch of the reduced case. See Fig.~\ref{fig: sketch full} for details. The sampling line is composed of two parts, and is reflected at $x = L/2$. The asymmetrical variables ($u_y$, and $\sigma_{xy}$) are transformed by the following function: $u(x) = -u, \ x < L/2; = u, \ x \geq L/2$.}\label{fig: sketch half}
\end{figure}

The results of computations are analyzed along the red dashed line, which is also shown on the sketch in Fig.~\ref{fig: sketch full}.
As meshless discretization is used, the results along the line are interpolated using the Sheppard interpolation from Python's \verb*|photutils| package. To stabilize the computations, and to reduce computational time, the case is further simplified by reducing the domain size along the symmetry line ($u_y = 0$, $\frac{\partial u_x}{\partial y} = 0$), producing the problem setup to that shown in Fig.~\ref{fig: sketch half}. To stabilize the computations even further, the corner nodes at Neumann boundary conditions (traction boundaries) are removed.

The computations are performed at different discretization resolutions, and different amount of load steps to display the model's convergence. Details on discretization resolutions are presented in Table~\ref{tab: discretization}.
Results are computed for each discretization, using $N_{\text{load}} \in \left \{ 5, 10, 50, 100 \right \}$ load steps.

\begin{table}
	\caption{Different discretization densities $\frac{dx}{L}$, and corresponding number of generated nodes in a domain.}\label{tab: discretization}
	\begin{center}
		\renewcommand\arraystretch{2}
		\begin{tabular}{|c|c|}
			\hline
			\rule{0pt}{3ex} \textbf{$\frac{dx}{L} \ [-]$} & \textbf{No. nodes} \\
			\hline
			\rule{0pt}{3ex} $\frac{1}{19}$                & 110                \\
			\hline
			\rule{0pt}{3ex} $\frac{1}{49}$                & 659                \\
			\hline
			\rule{0pt}{3ex} $\frac{1}{99}$                & 2585               \\
			\hline
			\rule{0pt}{3ex} $\frac{1}{149}$               & 5745               \\
			\hline
			\rule{0pt}{3ex} $\frac{1}{199}$               & 10204              \\
			\hline
			\rule{0pt}{3ex} $\frac{1}{249}$               & 15876              \\
			\hline
			\rule{0pt}{3ex} $\frac{1}{299}$               & 22852              \\
			\hline
		\end{tabular}
	\end{center}
\end{table}

The entire implementation is written in C++, using our in-house developed MEDUSA library~\cite{slak2021}. MEDUSA's built-in fill and relax algorithms~\cite{slak2019} are used to create the uniform irregular computational domains.
Differential operators are approximated with the RBF-FD and computed using the support of $n=50$ nearest points.
Eq. (\ref{eq: nce}) is solved implicitly with the \verb*|BiCGSTAB| solver from the Eigen library for linear algebra~\cite{eigenweb}.
In the main loop, the implicit solver operates sequentially, but local operations, such as computation of stress tensor $\sigma$, and the local iterations are executed in parallel with OpenMP.
The compiling is done with the \verb*|g++| compiler, using the following flags: \verb*|-Wall| \verb*|-O3| \verb*|-march=native| \verb*|-fopenmp| \verb*|-std=c++17|.
The results were exported to text-VTK files, and processed with Python and ParaView.
Computations were performed on a laptop, with an Intel Core i7-8750H CPU, and 16 GB DDR4 RAM.

For comparison and validation purposes, the same problem is also solved using the commercial software Abaqus FEA, which utilizes FEM. Abaqus solution is obtained with the full domain setup (Fig.~\ref{fig: sketch full}) using $N_{\text{FE}}=5000$ 8-node biquadratic plane stress quadrilateral (\verb*|CPS8R|) finite elements resulting in total of 15301 nodes. Finite element type \verb*|CPS8R| was chosen to obtain a finite element with 4 integration points and thus reducing the error of mapping the secondary variables to nodal values and additionally avoiding the potential Hourglass and shear locking issues. Reduced integration mode additionally reduces the computational times.

\section{Results}
\begin{figure}[htbp]
	\centering
	\includegraphics[width=\columnwidth]{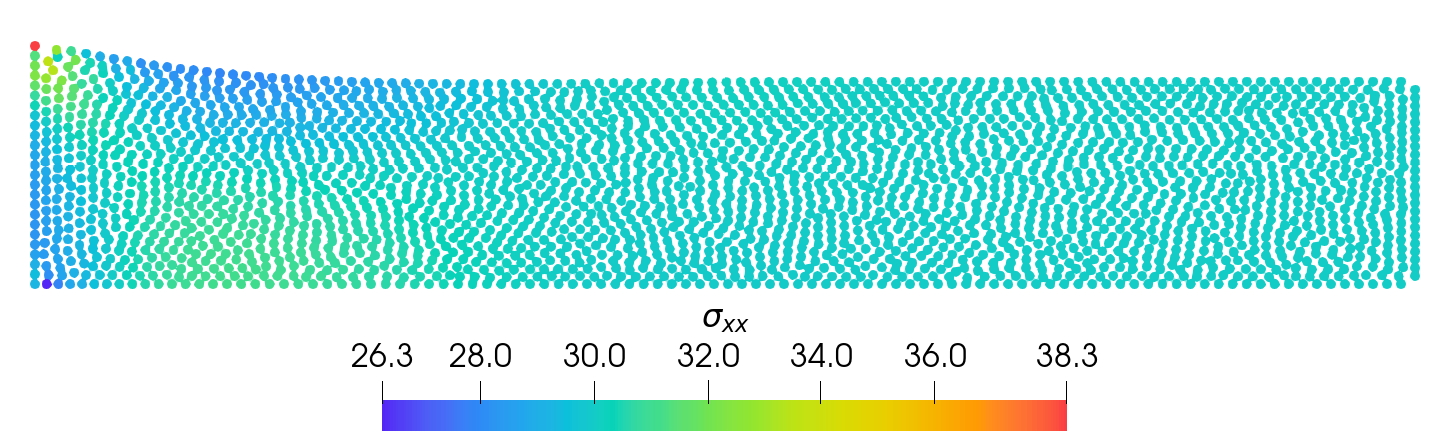}
	\caption{Visualization of the solution in a domain with 2585 nodes after 10/10 load steps. The nodes are displaced by $50 \bm{u}$, and the color-map is showing $\sigma_{xx}$.}\label{fig: solution}
\end{figure}

\begin{figure}[htbp]
	\centering
	\includegraphics[width=\columnwidth]{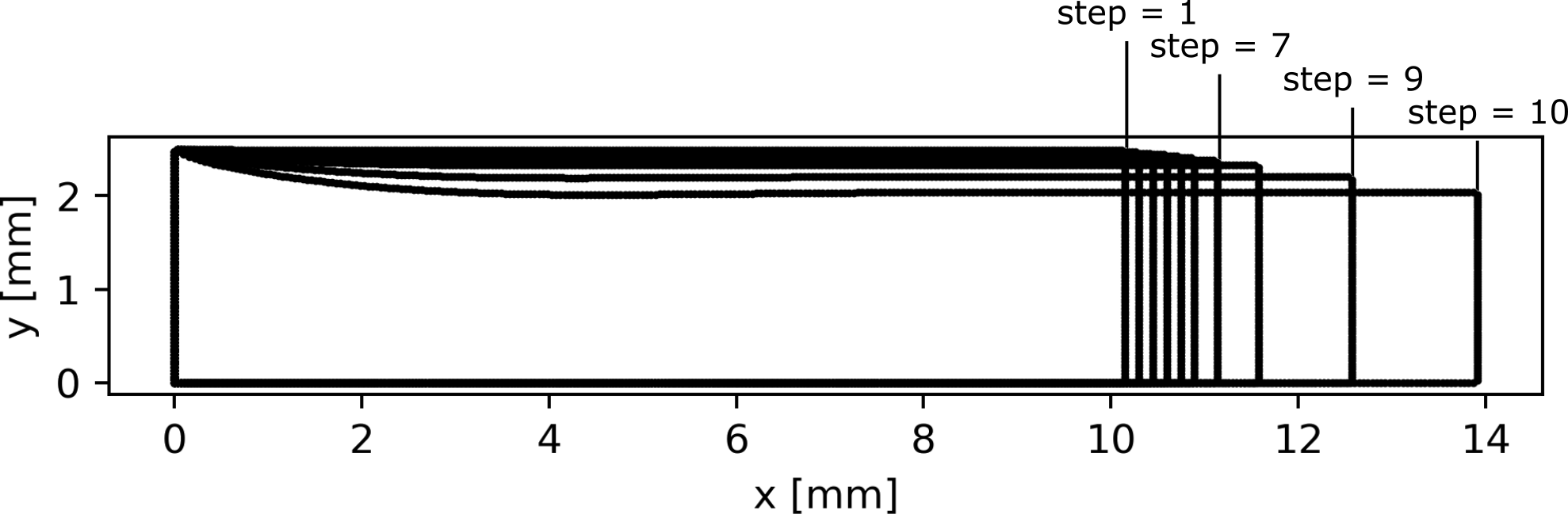}
	\caption{Visualization of the solution in a domain with 22852 nodes, showing the boundary nodes positions, displaced by $50 \bm{u}$, at load steps 1/10 - 10/10.}\label{fig: stretch}
\end{figure}

Example results of computations are shown in Figs.~\ref{fig: solution} and~\ref{fig: stretch}.
Fig.~\ref{fig: solution} is displaying the $\sigma_{xx}$ field on the discretized domain, with nodes displaced by $50 \bm{u}$.
A higher stress concentration can be observed in the north-western corner of the domain, which is surrounded by areas of lower $\sigma_{xx}$ stress, which then homogenizes along the body toward the eastern edge.
In Fig.~\ref{fig: stretch} it can be observed, how the material begins to yield, as it deforms plastically -- more drastic deformation increases can be observed after first couple of linear-elastic steps.

\begin{figure}[htbp]
	\centering
	\begin{subfigure}{0.49\textwidth}
		\centering
		\includegraphics[width=\textwidth]{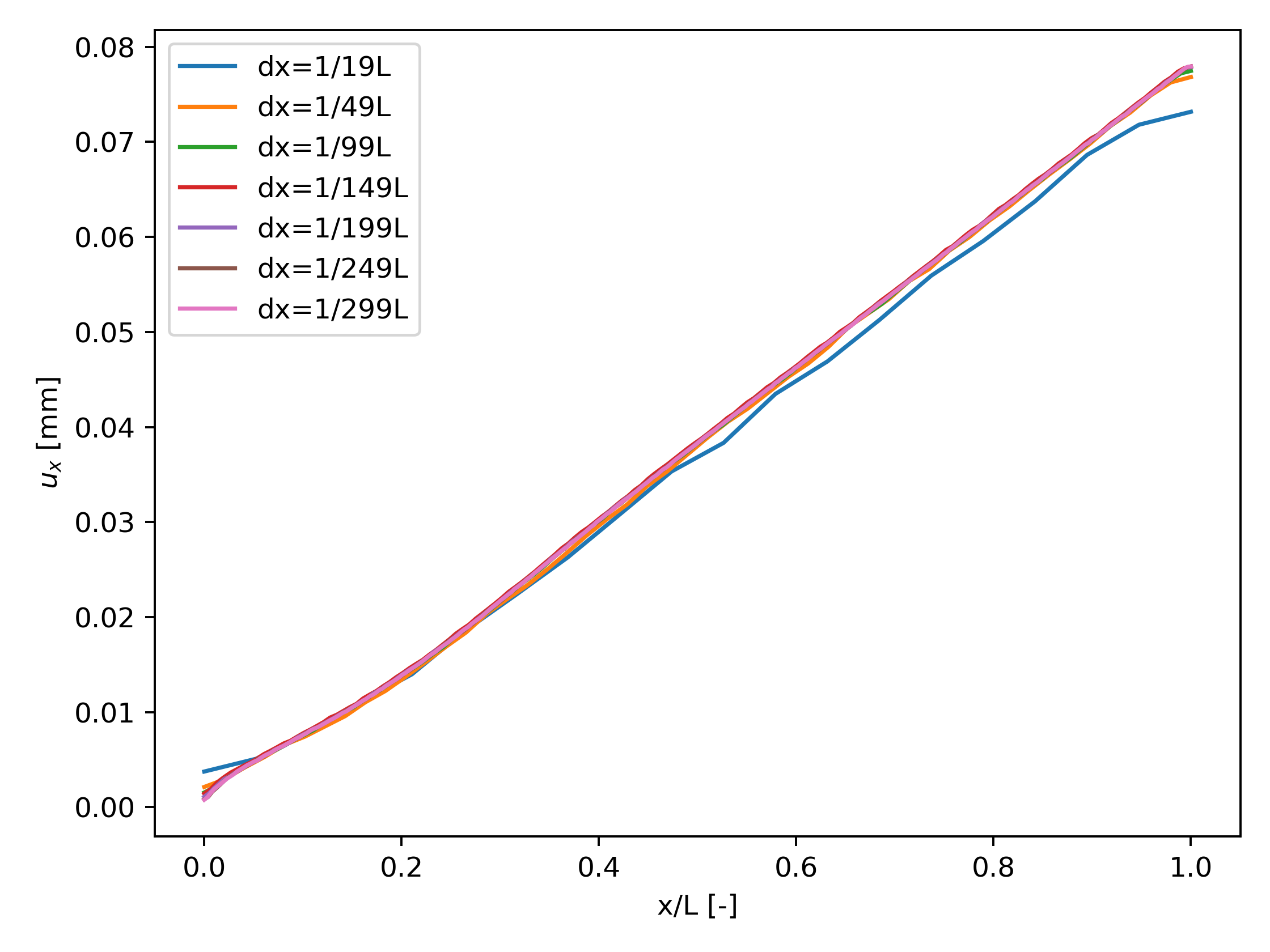}
		\caption{$u_x$}\label{fig: convdxa}
	\end{subfigure}
	\hfill
	\begin{subfigure}{0.49\textwidth}
		\centering
		\includegraphics[width=\textwidth]{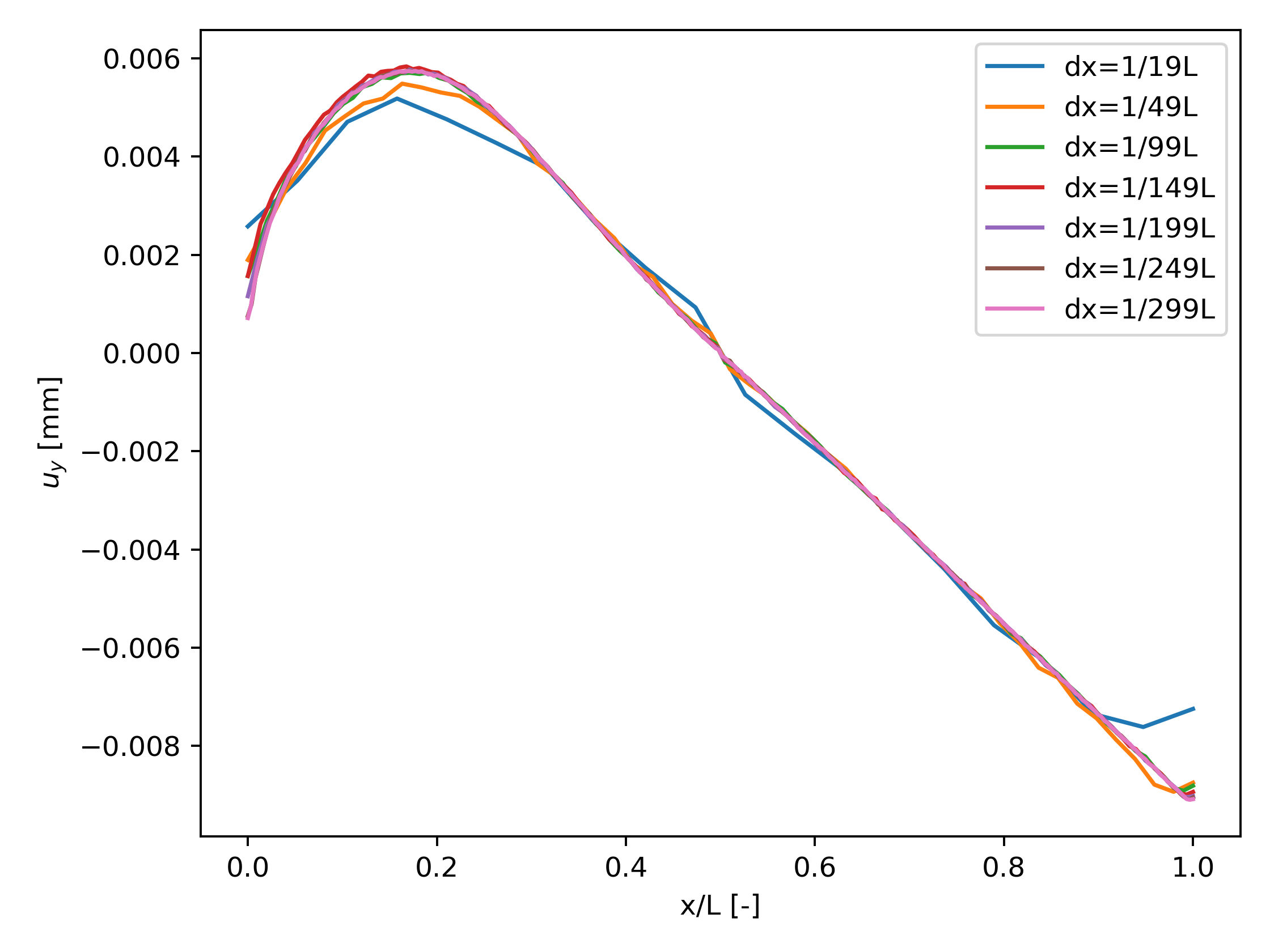}
		\caption{$u_y$}\label{fig: convdxb}
	\end{subfigure}
	\caption{Solution convergence due to varying discretization density. $u_x$ and $u_y$ plotted against dimensionless coordinate $x/L$, $N_{\text{steps}} = 10$.}\label{fig: convdx}
\end{figure}

\begin{figure}[htbp]
	\centering
	\begin{subfigure}{0.49\textwidth}
		\centering
		\includegraphics[width=\textwidth]{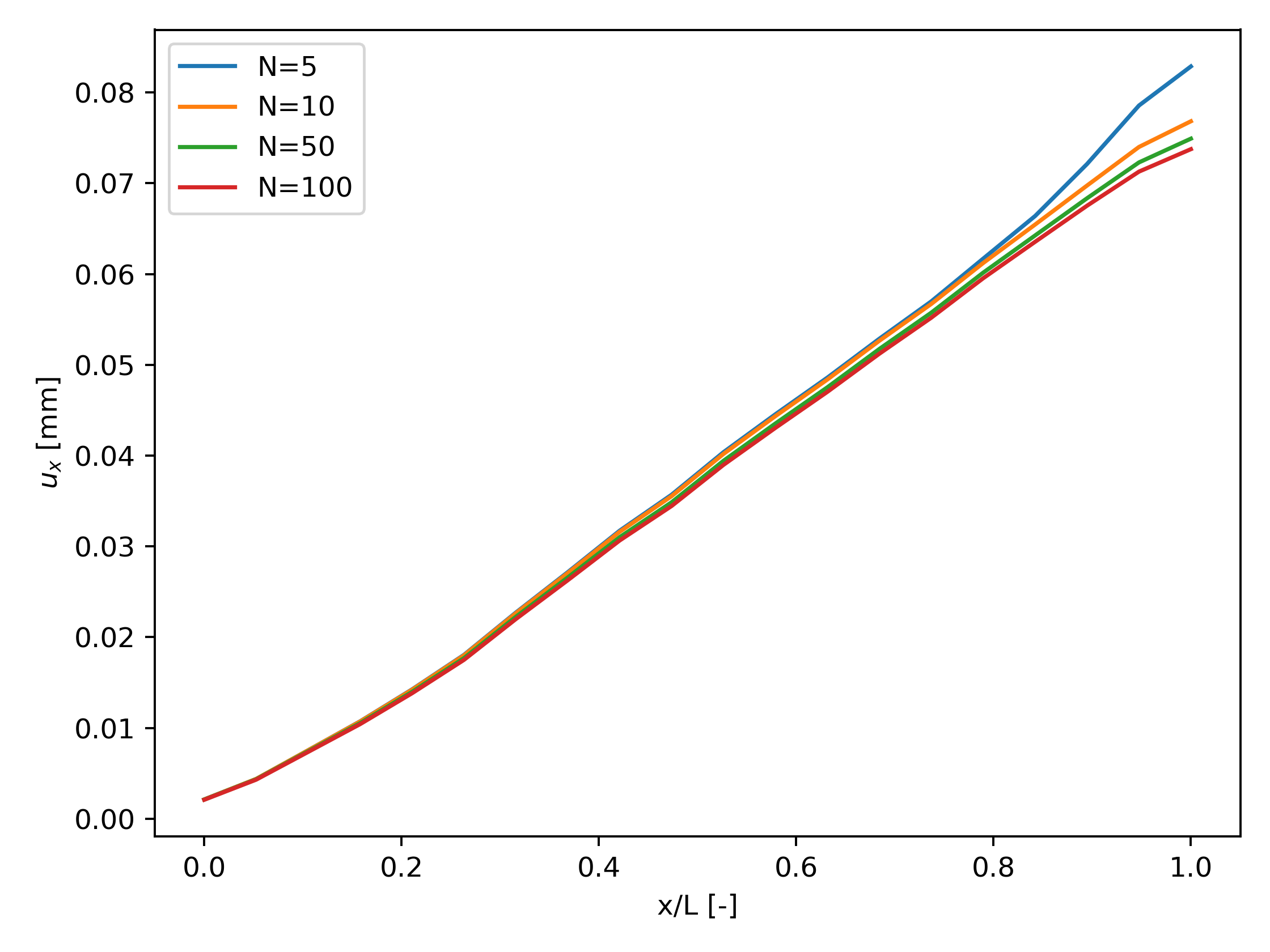}
		\caption{$u_x$, 259 nodes.}\label{fig: convna}
	\end{subfigure}
	\hfill
	\begin{subfigure}{0.49\textwidth}
		\centering
		\includegraphics[width=\textwidth]{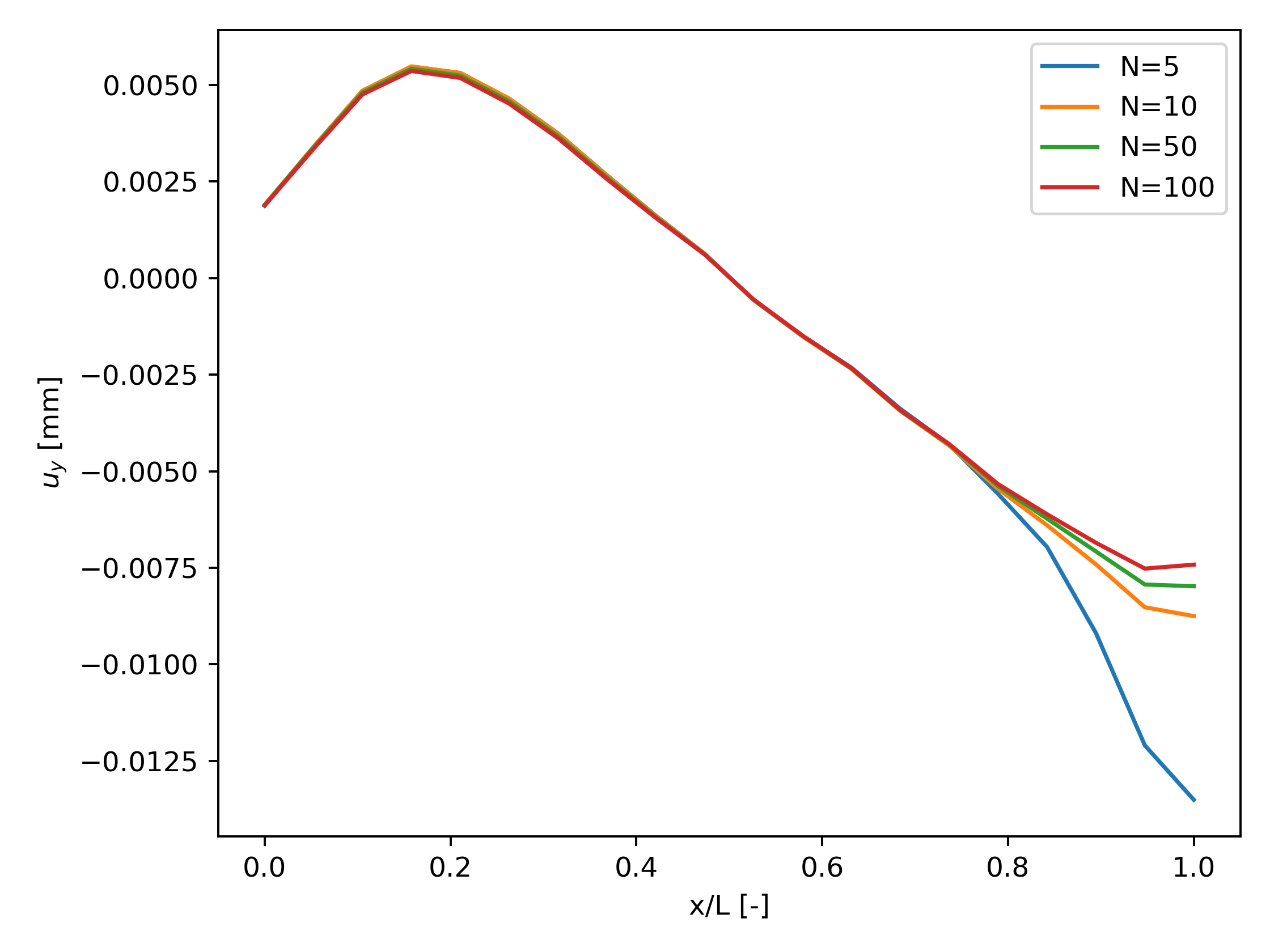}
		\caption{$u_y$, 259 nodes}\label{fig: convnb}
	\end{subfigure}
	\caption{Solution convergence due to varying number of load steps $N$ density. $u_x$ and $u_y$ plotted against dimensionless coordinate $x/L$, results for 259 nodes.}\label{fig: convn}
\end{figure}

The model's convergence by increasing the node density can be seen in Fig.~\ref{fig: convdx}, and in Fig.~\ref{fig: convn} is showing the load-step dependent convergence for $\bm{u}$ at low (659 nodes) discretization density.

\begin{figure}[htbp]
	\centering
	\begin{subfigure}{0.49\textwidth}
		\centering
		\includegraphics[width=\textwidth]{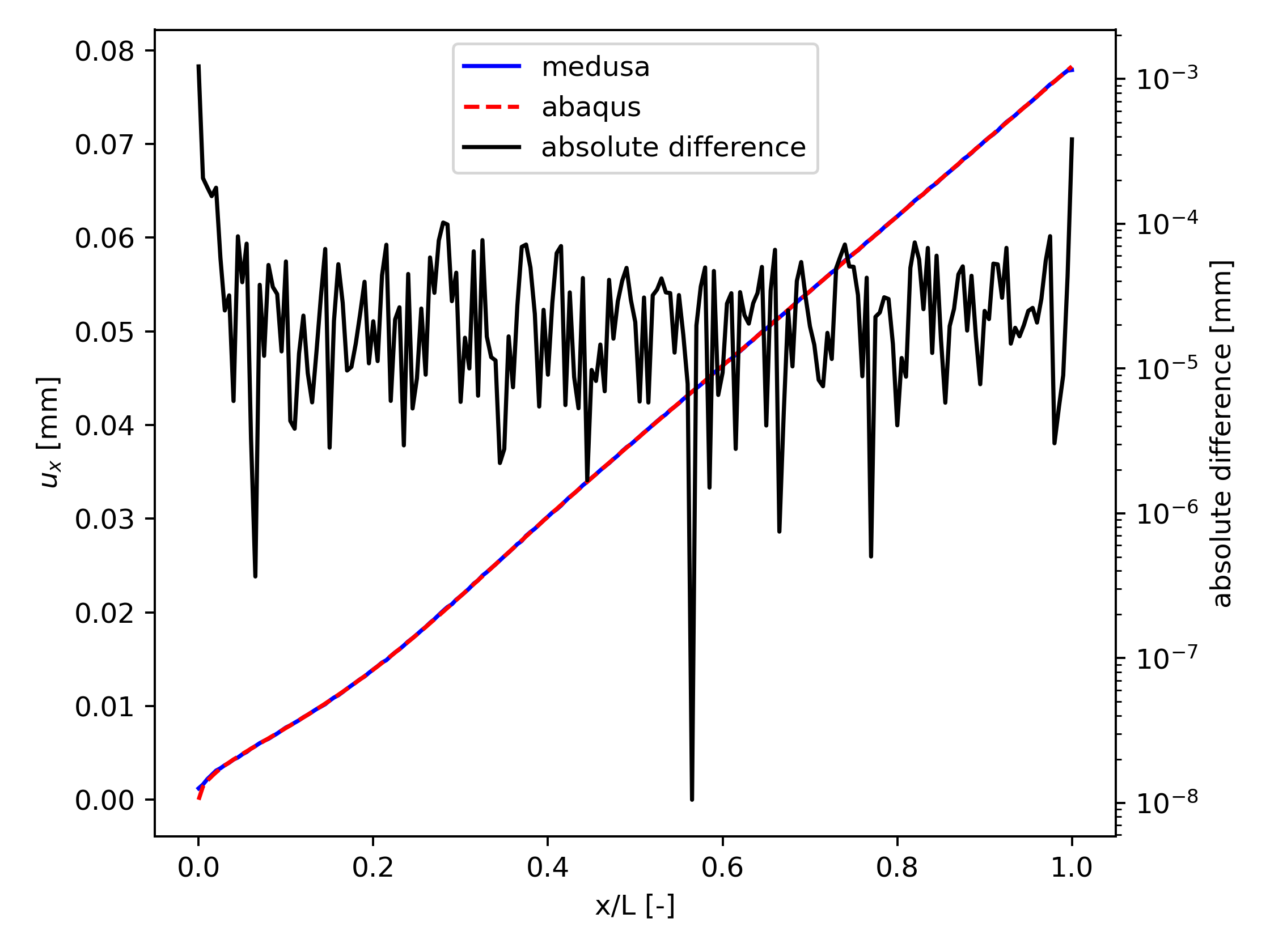}
		\caption{$u_x$}\label{fig: compua}
	\end{subfigure}
	\hfill
	\begin{subfigure}{0.49\textwidth}
		\centering
		\includegraphics[width=\textwidth]{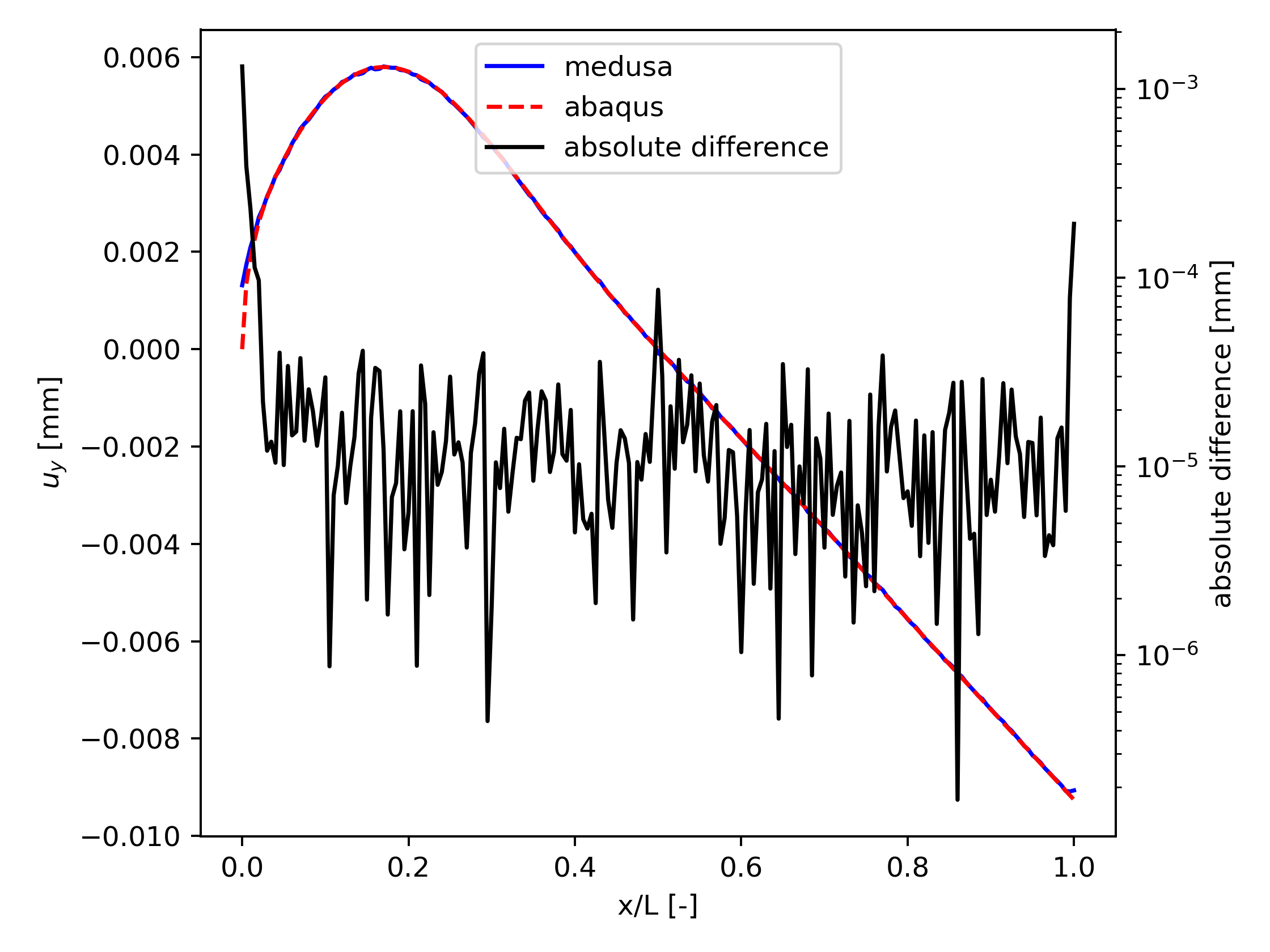}
		\caption{$u_y$}\label{fig: compub}
	\end{subfigure}
	\caption{Comparison of both components of $\bm{u}$, and the absolute difference between Medusa and Abaqus results plotted against dimensionless coordinate $x/L$.}\label{fig: compu}
\end{figure}

\begin{figure}[htbp]
	\centering
	\begin{subfigure}{0.49\textwidth}
		\centering
		\includegraphics[width=\textwidth]{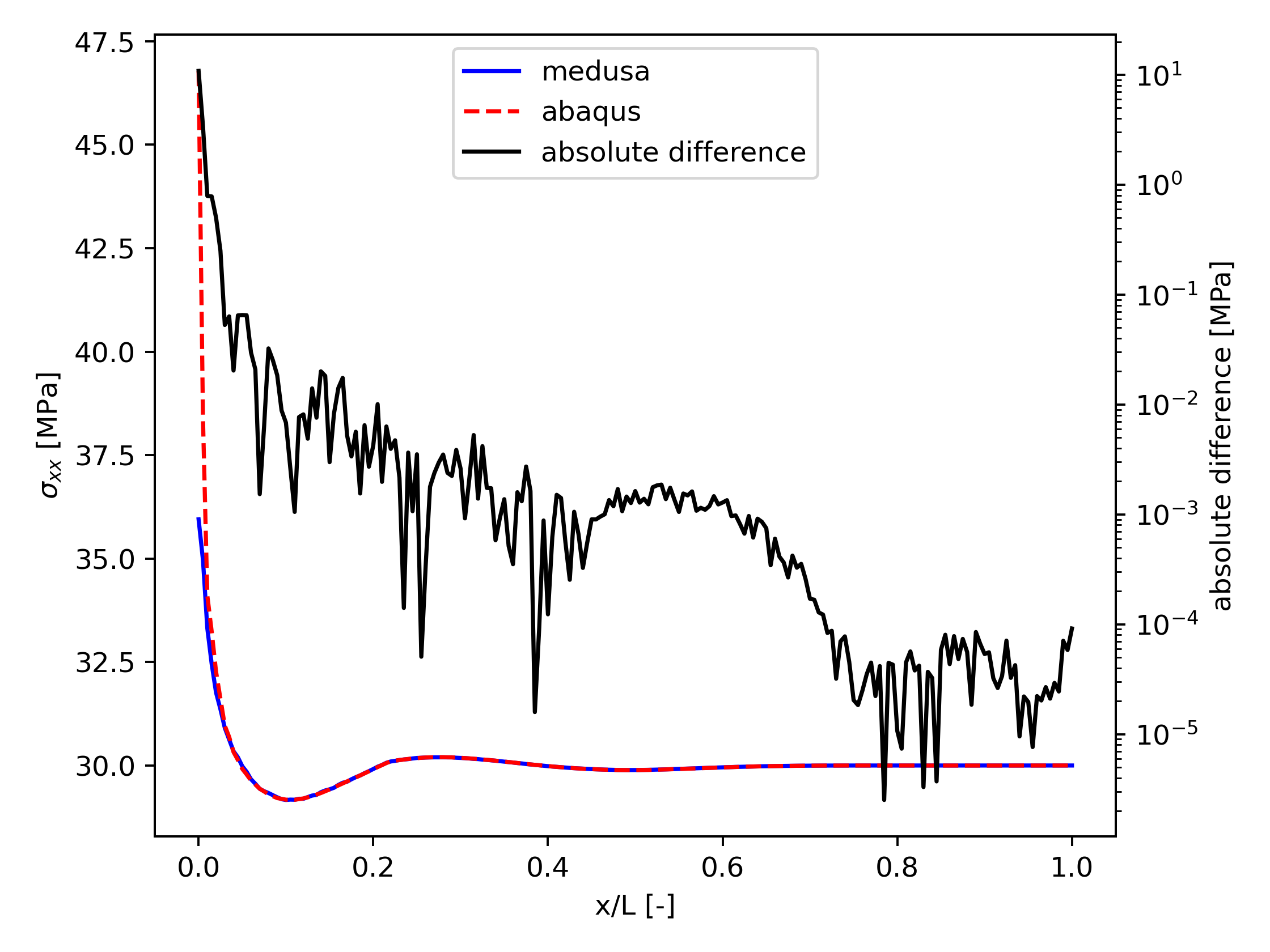}
		\caption{$\sigma_{xx}$}\label{fig: compsa}
	\end{subfigure}
	\hfill
	\begin{subfigure}{0.49\textwidth}
		\centering
		\includegraphics[width=\textwidth]{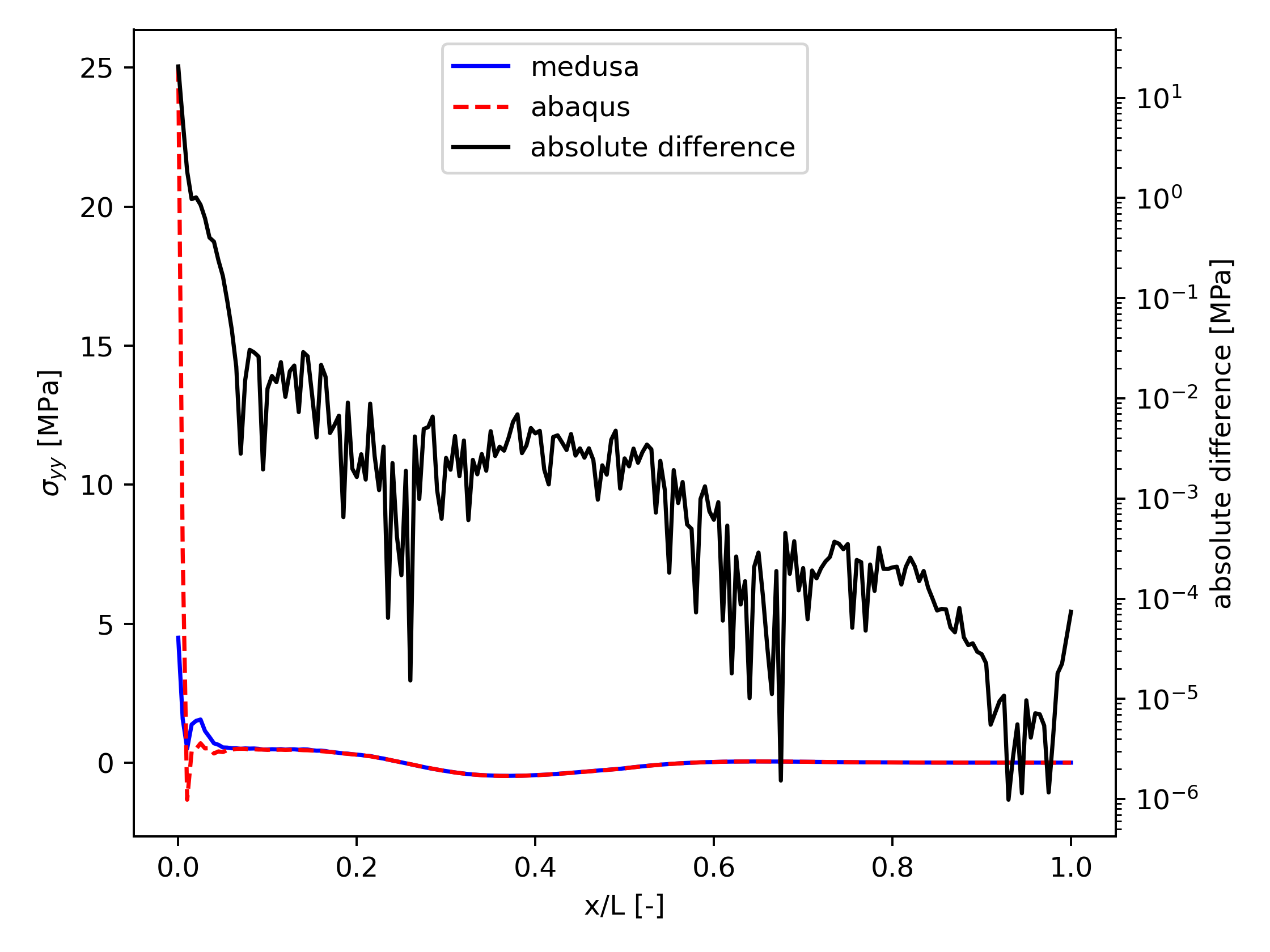}
		\caption{$\sigma_{yy}$}\label{fig: compsb}
	\end{subfigure}
	\hfill
	\begin{subfigure}{0.49\textwidth}
		\centering
		\includegraphics[width=\textwidth]{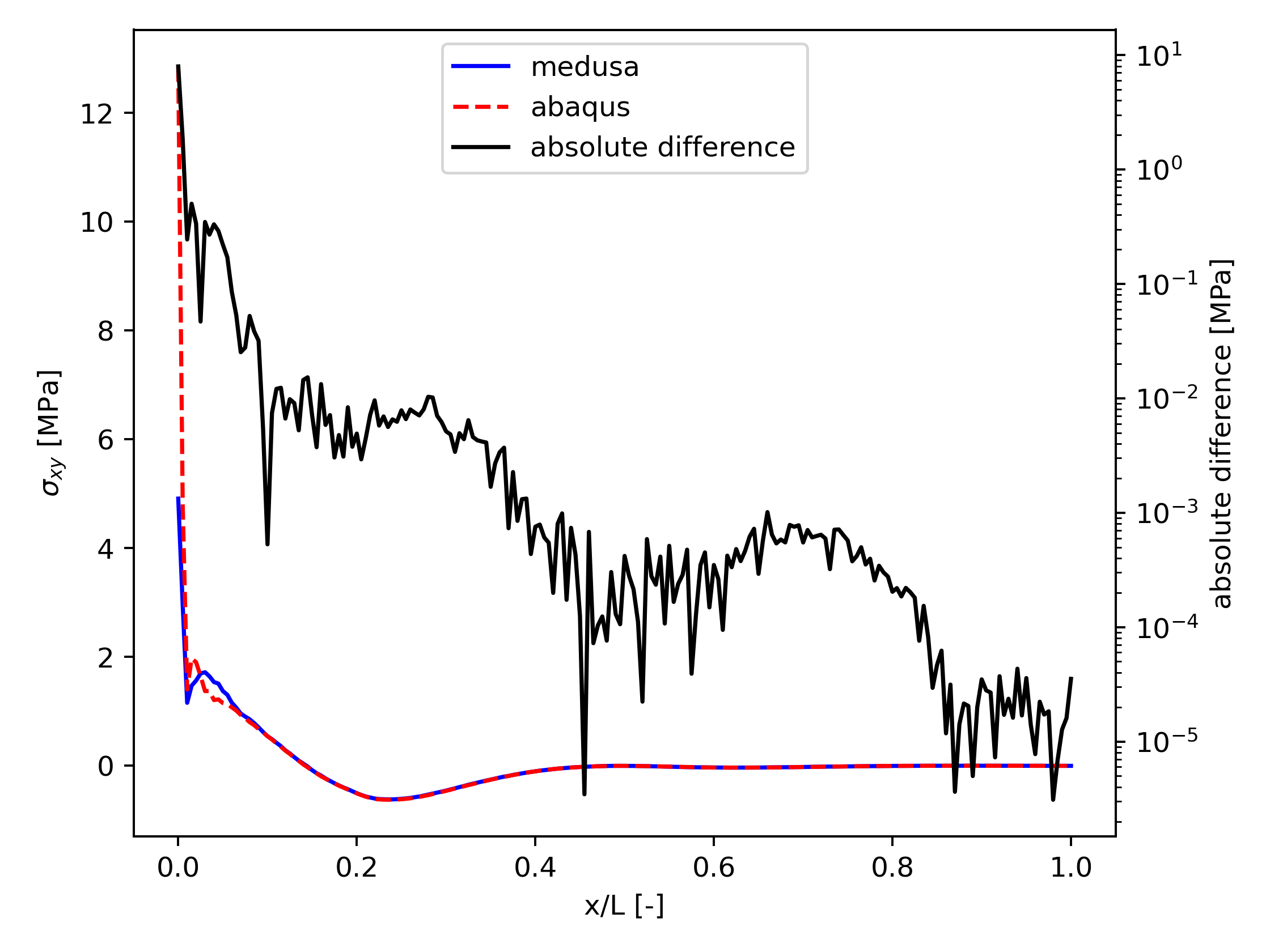}
		\caption{$\sigma_{xy}$}\label{fig: compsc}
	\end{subfigure}
	\caption{Comparison of the three components of $\sigma$, and the absolute difference between Medusa and Abaqus results plotted against dimensionless coordinate $x/L$.}\label{fig: comps}
\end{figure}

The solution with $N=10204$ computational nodes and $N_{\text{steps}} = 100$ load steps is compared with the Abaqus FEA solution along the aforementioned diagonal. The comparison is made for solutions of $\bm{u}$ (Fig.~\ref{fig: compu}), and $\sigma$ (Fig.~\ref{fig: comps}). In both cases the absolute difference plotted on the log scale represents $\Delta u = \left|u_{medusa} - u_{abaqus}\right|$, where indices ``medusa'' and ``abaqus'' are representing the meshless and Abaqus solutions, respectively.

\begin{figure}[htbp]
	\centering
	\includegraphics[width=\columnwidth]{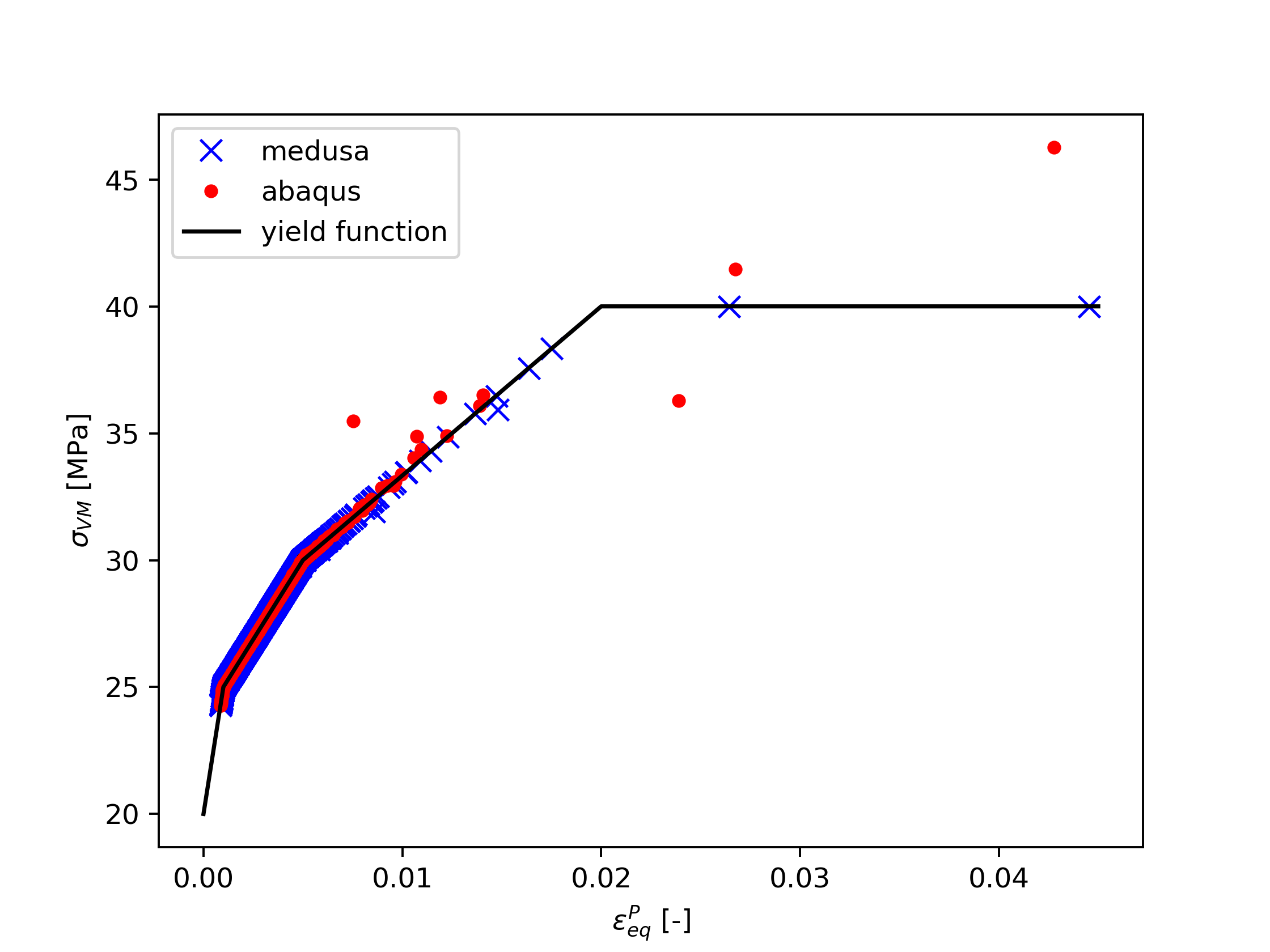}
	\caption{$\sigma_{VM}$ plotted against $\varepsilon^P_{eq}$ for meshless and Abaqus solutions. The solid line represents the yield function, which is flat beyond the (0.02, \SI{40}{\mega\pascal}) point, resulting in perfect plastic behaviour.}\label{fig: yieldf}
\end{figure}

Finally, all points from meshless (the same case as previously) and Abaqus solutions are plotted on the $\sigma_{VM}$ against $\varepsilon^P_{eq}$ graph along with the yield function in Fig.~\ref{fig: yieldf}, to analyze how well they follow the yield function.

The convergence plots in Fig.~\ref{fig: convn} show that the quality of the solution may increase by increasing the number of load steps, but at higher density discretizations this may be insignificant.
Finding the balance between discretization density and number of load steps for optimal accuracy and computational time is beyond the scope of this paper.

The comparison between the presented implementation of von Mises plasticity model with the commercial solver Abaqus (Figs.~\ref{fig: compu} and~\ref{fig: comps}) shows good agreement between both solutions.
This is shown by closely matched solution for $\bm{u}$ as is shown in Fig.~\ref{fig: compu}, and a similar $\sigma$ prediction in Fig.~\ref{fig: comps}.
One should note however, that the meshless domain had the corner nodes removed to stabilize the computation, which resulted in greater difference between the two solutions especially at 0 $x$-coordinate.
The difference between the two solutions for $\sigma$ decreases as one is going up the $x$-coordinate, except for the corner node, for reasons mentioned above.
The decrease in the difference is likely due to $\sigma_{xx}$ being prescribed at $x = L$ through the traction boundary condition, and both approaches fulfill it well.
A source of error, which may contribute to the noisy nature of the absolute difference plots, is the interpolation, which was used in post processing to extract the data along the prescribed lines.

Fig.~\ref{fig: yieldf} is showing that the implemented plasticity model closely follows the yield function.
However, some of the values obtained by both approaches appear to be violating the yield criterion.
In the meshless case these points are located on the fixed edge, and in Abaqus' case these are at both singularities at the fixed edge or in their close proximity.

\section{Conclusion}
This work demonstrates the use of RBF-FD in modeling small-strain plasticity.
Under the plane stress assumption, a simple 2D problem is solved, using the von Mises plasticity model.
For comparison, an identical problem is also solved in Abaqus, which uses FEM.
Results show good agreement between both approaches.
The analyzed example case also indicates that increasing the number of load steps may stabilize the solution in a scarcely populated domain, whereas this effect is not as pronounced in a densely populated one.
Generally it can be concluded that RBF-FD can substitute FEM for solving similar small-strain plasticity problems.

\bibliographystyle{IEEEtran}
\bibliography{refs}

\end{document}